\documentclass[12pt,a4]{article}
\usepackage{indentfirst}
\usepackage{amssymb}
\newtheorem{theo}{THEOREM}
\newtheorem{p}[theo]{PROPOSITION}
\newcommand{\h}{Theorem }
\newcommand{\hh}{[C] Theorem 5.6.}
\newcommand{\pr}{Proposition }
\newcommand{\prr}{[C] Proposition 5.6.}
\newcommand{\qedd}{\eqno{\rule{3mm}{3mm}}}
\newcommand{\qed}{\hfill {\rule{3mm}{3mm}}}
\newcommand{\ma}[5]{\ensuremath{#1:#2\longrightarrow #3\,, \quad #4 \longmapsto  #5}}
\newcommand{\mad}[4]{\ensuremath{#1\longrightarrow #2\,, \quad #3 \longmapsto  #4}}
\newcommand{\mac}[3]{\ensuremath{#1:#2\longrightarrow #3}}

\newcommand{\ri}[1]{Hilbert right $#1$-module}
\newcommand{\ris}[1]{selfdual\, Hilbert\, right\, $#1$-module}
\newcommand{\f}{\ensuremath{\mathcal{F}}}

\newcommand{\me}[2]{\ensuremath{\left\{\,\left. #1\;\right|\; #2 \, \right\}}}
\newcommand{\s}[2]{\ensuremath{\left \langle  \, #1\left|\, #2 \right. \, \right \rangle}}
\newcommand{\sa}[2]{\ensuremath{\left \langle \, #1 \,, \, #2 \, \right \rangle}}
\newcommand{\si}[1]{\ensuremath{\sum\limits_{#1}}}
\newcommand{\cb}[1]{\ensuremath{\mathop{\bigcirc \hspace{-2.7mm} |}\limits_{ \hspace{2mm} #1}}}
\newcommand{\lb}[2]{\ensuremath{\mathcal L_{#1}(#2)}}
\newcommand{\lc}[3]{\ensuremath{\mathcal L_{#1}(#2,#3)}}
\newcommand{\ld}[3]{\ensuremath{\hat \mathcal L_{#1}(#2,#3)}}

\newcommand{\bk}{\ensuremath{\mathrm{I\! K}}}
\newcommand{\aaaa}[4]{\ensuremath{ \left\{ \begin{array}{c@{\quad \mbox{if} \quad}c} #1&#2 \\
#3&#4 \end{array}} \right.}
\newcommand{\ti}[3]{\ensuremath{\widetilde{\overbrace{(#1,#2,#3)}}}}
\newcommand{\tit}[2]{\ensuremath{\widetilde{(#1,#2)}}}
\renewcommand{\labelenumi}{\alph{enumi})}

\begin{document}

\begin{center}
\Large {\bfseries
THE DOUBLE COMMUTATION THEOREM FOR SELFDUAL HILBERT RIGHT W*-MODULES} 
\end{center}

\begin{center}
CORNELIU CONSTANTINESCU
\end{center}
\vspace{1cm}
\begin{abstract}
Let $E$ be a W*-algebra, $H$ a selfdual Hilbert right $E$-module, and $\f$ an involutive unital subalgebra of $\lb{E}{H}$. We prove that the double commutant of $\f$ is the W*-subalgebra of $\lb{E}{H}$ generated by $\f$ (\h\ref{8063}). The proofs work simultaneously for the real and for  the complex case. If $E$ is the field of real or complex numbers then $H$ is a Hilbert space and this result becomes the well-known classical double commutation theorem. 
\end{abstract}
{\it
AMS classification code: 46L08, 46L10 
 
\hspace{-26pt} Keywords: \ris{W^*}s, the double commutation theorem}

In general we use the notation and terminology of [C]. In the sequel we give a list of some notation used in this paper.

\renewcommand{\labelenumi}{\arabic{enumi}.}
\begin{enumerate}
\item  $\bk$ denotes the field of real or the field of complex numbers. The whole theory is developed in parallel for the real and complex case (but the proofs coincide).
\item  For all $i,j$ we denote by $ \delta _{ij}$ the Kronecker's symbol:
$$\delta _{ij}:= \aaaa{1}{i=j}{0}{i\not=j}\,.$$
\item If $E,F$ are vector spaces in duality then $E_F$ denotes the vector space $E$ endowed with the locally convex topology of pointwise convergence on $F$, i.e. with the weak topology $\sigma (E,F)$.
\item If $E$ is a Banach space then $E'$ denotes its dual and $E^{\#}$ its unit ball :
$$ E^{\#}:= \me{x\in E}{\|  x\|  \leq 1}.$$
\item If $E,F$ are Banach spaces then $\lc{}{E}{F}$ denotes the Banach space of operators from $E$ to $F$. 
\item If $E$ is a C*-algebra then $Pr\,E$ denotes the set of orthogonal projections of $E$. If $E$ is a W*-algebra then $\ddot{E} $ denotes its predual.
\item Let $E$ be a W*-algebra. 
A C*-subalgebra $F$ of $E$ is called a W*-subalgebra of $E$ if it is a W*-algebra and if the inclusion map $F_{\ddot F}\rightarrow E_{\ddot E}$ is continuous. We say that an involutive subalgebra $F$ of $E$ generates $E$ as a W*-subalgebra if every W*-subalgebra of $E$ containing $F$ is equal to $E$; by [C] Corollary 4.4.4.12 a), this is equivalent to the assertion that $F$ is dense in $E_{\ddot E}$.
\item Let $E$ be a C*-algebra and $H,K$ \ri{E}s.  $\lc{E}{H}{K}$ denotes the Banach subspace of $\lc{}{H}{K}$ of adjointable operators, $\ld{E}{H}{K}$ the Banach space of  $u\in \lc{}{H}{K}$ such that
$$(\xi ,x)\in H\times E\,\Longrightarrow\, u(\xi x)=(u\xi )x,$$ 
$$\lb{E}{H}:=\lc{E}{H}{H},\qquad \hat H:=\ld{E}{H}{E}$$
and
$$\ma{\s{\cdot}{\xi }}{H}{E}{\eta }{\s{\eta }{\xi }}$$ 
for every $\xi \in H$. If for every $u\in \hat H$ there is a $\xi \in H$ with
$$u:=\s{\cdot}{\xi }$$
then $H$ is called selfdual. 
\item Let $E$ be a W*-algebra and $H$ a \ri{E}. We put for $a\in \ddot E$ and $\xi ,\eta \in H$,
$$\ma{\widetilde {(a,\xi )}}{H}{\bk}{\zeta }{\sa{\s{\zeta }{\xi }}{a}},$$
$$\ma{\ti {a}{\xi}{\eta}}{\lb{E}{H}}{\bk}{u}{\sa{\s{u\xi }{\eta }}{a}}$$
and denote by $\ddot H$ (by $\stackrel{...}{H}$) the closed vector subspace of $H'$ (of $\lb{E}{H}'$) generated by 
$$\me{\widetilde {(a,\xi )}}{(a,\xi )\in \ddot E\times H}  \; \left(\mbox{by} 
\me{\ti {a}{\xi}{\eta}}{(a,\xi ,\eta )\in \ddot E\times H\times H}\right)\;.$$
If $H$ is selfdual then $\ddot H$ (resp. $\stackrel{...}{H}$) is the predual of $H$ (resp. of $\lb{E}{H}$)([C] \pr 5.6.3.3, [C] \h 5.6.3.5 a)).
\item If $E$ is a C*-algebra, $H$ a \ri{E}, and $I$ a finite set then $\cb{i\in I}H$ denotes the \ri{E} obtained by endowing the vector space $H^I$ with the right multiplication
$$\mad{H^I\times E}{H}{(\xi ,x)}{(\xi _ix)_{i\in I}}$$
and with the inner product
$$\mad{H^I\times H^I}{E}{(\xi ,\eta )}{\si{i\in I}\s{\xi _i}{\eta _i}}.$$
\item $\f^c$ denotes the commutant and $\f^{cc}$ the double commutant of $\f$.
\end{enumerate}
 
 \renewcommand{\labelenumi}{\alph{enumi})}
 
\begin{p}\label{7445}
If $E$ is a C*-algebra, $H$ a \ris{E}, and $K$ a Hilbert right $E$-submodule of $H$, then the following are equivalent:
\begin{enumerate}
\item K is selfdual.
\item There is a $p\in Pr \,\lb{E}{H}$ with $K=p(H)$.

If $E$ is a W*-algebra then the above assertions are equivalent to the following one:
\item $K$ is closed in $H_{\ddot{H} }$.
\end{enumerate}
\end{p}
 
 $a\Rightarrow  b.$ Let $u$ be the inclusion map $K\rightarrow H$. Then
 $$u\in \ld{E}{K}{H}=\lc{E}{K}{H}  $$
 (\prr 2.4). For $(\xi ,\eta )\in H\times K$,
 $$\s{\xi }{\eta }=\s{\xi }{u\eta }=\s{u^*\xi }{\eta }$$
 so 
 $$\xi \in K\Longrightarrow u^*\xi =\xi $$
 and
 $$\s{u^*uu^*\xi }{\eta }=\s{uu^*\xi }{\eta }=\s{u^*\xi }{\eta },$$
 $$u^*uu^*\xi =u^*\xi ,\qquad u^*uu^*=u^*,\qquad uu^*uu^*=uu^*,$$
 $$p:=uu^*\in Pr\,\lb{E}{H},\qquad p(H)=K.$$
 
$b\Rightarrow a.$ Let $u\in \ld{E}{K}{E}$. Put \mbox{$v:=u\circ p\in \ld{E}{H}{E}$}. Since $H$ is selfdual there is an $\eta \in H$ with
$$v=\s{\cdot}{\eta }\,.$$
For $\xi \in K$,
$$u\xi =up^2\xi =vp\xi =\s{p\xi }{\eta }=\s{\xi }{p\eta }$$
so
$$u=\s{\cdot}{p\eta }$$
and $K$ is selfdual. 

$b\Rightarrow c.$ By \prr3.4 c), the map
$$\mac{p}{H_{\ddot{H} }}{H_{\ddot{H} }}$$
is continuous, so
$$K=\me{\xi \in H}{\xi =p\xi }$$
is closed in $H_{\ddot{H} }$.

$c\Rightarrow a.$ Since $H^{\#}_{\ddot{H} }$ is compact ( \prr 3.3 $a\Rightarrow b$), it follows that $K_{\ddot{K}}^{\#}=K_{\ddot{H}}^{\#} $ is also compact, so $K$ is selfdual (\prr3.3 $b\Rightarrow a$).\qed

\begin{p}\label{8060}
Let $E$ be a W*-algebra, $H$ a \ris{E}, $\f$ an involutive subalgebra of $\lb{E}{H}$, $\xi \in H$, and $K$ the closure in $H_{\ddot{H} }$ of $\me{u\xi }{u\in \f}$. Then $w(K)\subset K$ for every $w\in \f^{cc}$.
\end{p}

Let $\xi _0\in K$ and $v\in \f$. Let $(a_i,\eta _i)_{i\in I}$ be a finite family in $\ddot{E}\times H $. There is a $u\in \f$ such that
$$\left|{\sa{u\xi -\xi _0}{\tit{a_i}{v^*\eta _i}}}\right|<1.$$
for every $i\in I$. Then
$$\left|\sa{vu\xi -v\xi _0}{\tit{a_i}{\eta _i}}\right|=|\sa{\s{vu\xi -v\xi _0}{\eta _i}}{a_i}|=$$
$$=|\sa{\s{u\xi -\xi _0}{v^*\eta _i}}{a_i}|=\left|{\sa{u\xi -\xi _0}{\tit{a_i}{v^*\eta _i}}}\right|<1.$$
for every $i\in I$, so $v\xi _0\in K$. Thus $v(K)\subset K$. If $\eta \in K^\perp $ then for $\zeta \in K$,
$$\s{v\eta }{\zeta }=\s{\eta }{v^*\zeta }=0,$$
so $v\eta \in K^\perp $ and $v(K^\perp )\subset K^\perp $.

By \pr\ref{7445} c$\Rightarrow $b, there is a $p\in Pr\,\lb{E}{H}$ with $K=p(H)$. By the above, for $\eta \in H$,
$$pv\eta =pvp\eta +pv(\eta -p\eta )=pvp\eta =vp\eta ,$$
so $pv=vp$, $p\in \f^c$. Thus 
$$w\xi _0=wp\xi _0=pw\xi _0\in K,\qquad w(K)\subset K.\qedd$$

\begin{p}\label{8062}
Let $E$ be a C*-algebra, $H$ a \ri{E}, and $I$ a finite set. We put for every $u\in \lb{E}{H}$,
$$\ma{\tilde{u} }{\cb{i\in I}H}{\cb{i\in I}H}{\xi }{(u\xi_i )_{i\in I}}$$
and use the matrix notation of \emph{\prr4.16 d)}.
\begin{enumerate}
\item For every $u\in \lb{E}{H}$ and $i,j\in I$,
$$(\tilde{u}_{i,j} )=\delta _{i,j}u.$$
\item The following are equivalent for all $u\in \lb{E}{H}$ and $v\in \lb{E}{\cb{i\in I}H}$:
\begin{enumerate}
\item $\tilde{u}v=v\tilde{u} $. 
\item $uv_{i,j}=v_{i,j}u$ for all $i,j\in I$.
\end{enumerate}
\end{enumerate}
\end{p}

a) Let $\xi \in H$ and for every $j\in I$ put $\eta_j :=(\delta _{i,j}\xi )_{i\in I}$. Then for $i,j\in I$,
$$\delta _{i,j}u\xi =u(\delta _{i,j}\xi )=u((\eta_j) _i)=(\tilde{u}\eta_j  )_i=\si{k\in I}\tilde{u}_{i,k}(\eta_j) _k=\si{k\in I}\tilde{u}_{i,k}\delta _{k,j}\xi =\tilde{u}_{i,j}\xi ,   $$
so $\tilde{u}_{i,j}=\delta _{i,j}u $.

b) By a), for $i,j\in I$,
$$(\tilde{u}v )_{i,j}=\si{k\in I}\tilde{u}_{i,k}v_{k,j}=\si{k\in I}\delta _{i,k}uv_{k,j}=uv_{i,j}, $$
$$(v\tilde{u} )_{i,j}=\si{k\in I}v_{i,k}\tilde{u}_{k,j}=\si{k\in I}v_{i,k}\delta _{k,j}u=v_{i,j}u, $$
which proves the assertion.\qed

\begin{theo}\label{8063}
Let $E$ be a W*-algebra, $H$ a $\ris{E}$, and $\f$ an involutive unital subalgebra of $\lb{E}{H}$. Then $\f^{cc}$ is the W*-subalgebra of $\lb{E}{H}$ \emph{(\hh3.5 b))} generated by $\f$. 
\end{theo}

Let $u\in \f^{cc}$ and let $(a_i,\xi _i,\eta _i)_{i\in I}$ be a finite family in $\ddot{E}\times H\times H $. We put $K:=\cb{i\in I}H, $
$$\ma{\tilde{v} }{K}{K}{\zeta }{(v\zeta_i )_{i\in I}}$$
for every $v\in \lb{E}{H}$, and
$$\tilde{\f}:=\me{\tilde{v} }{v\in \f} .$$
Let $w\in (\tilde{\f})^c $. Let $v\in \f$ and $i,j\in I$. By \pr\ref{8062} $b_1\Rightarrow b_2$, $vw_{i,j}=w_{i,j}v$, so $w_{i,j}\in \f^c$. It follows $uw_{i,j}=w_{i,j}u$, so by \pr\ref{8062} $b_2\Rightarrow b_1$, $\tilde{u}w=w\tilde{u}  $. Thus $\tilde{u}\in (\tilde{\f})^{cc} $.

Put $\xi :=(\xi _i)_{i\in I}$, $\zeta _j:=(\delta _{i,j}\eta _i)_{i\in I}$ for every $j\in I$, and denote by $L$ the closure of $\me{\tilde{v}\xi  }{v\in \f}$ in $K_{\ddot{K} }$. Since $\xi \in L$ and $\tilde{u}\in (\tilde{\f} )^{cc} $ it follows from \pr\ref{8060}, that there is a $v\in \f$ with
$$\left|\sa{\tilde{u}-\tilde{v}  }{\ti{a_i}{\xi}{\zeta _i}}\right|<1$$
for every $i\in I$. For $j\in I$,
$$\s{(\tilde{u}-\tilde{v}  )\xi }{\zeta _j}=\si{i\in I}\s{(u-v)\xi _i}{(\zeta _j)_i}=$$
$$=\si{i\in I}\s{(u-v)\xi _i}{\delta _{i,j}\eta _i}=\s{(u-v)\xi _j}{\eta _j},$$
so for $i\in I$,
$$\left|\sa{u-v}{\ti{a_i}{\xi _i}{\eta _i}}\right|=|\sa{\s{(u-v)\xi _i}{\eta _i}}{a_i}|=$$
$$=\mid \sa{\s{(\tilde{u}-\tilde{v}  )\xi }{\zeta _i}}{a_i}\mid =\left|\sa{\tilde{u}-\tilde{v}  }{\ti{a_i}{\xi }{\zeta _i}}\right|<1.$$
Thus $u$ belongs to the closure $\mathcal{G}$ of $\f$ in $\lb{E}{H}_{\stackrel{...}{H}}$, so $\f^{cc}\subset \mathcal{G}$. Since $\f^{cc}$ is obviously closed in $\lb{E}{H}_{\stackrel{...}{H}}$ and contains $\f$, $\f^{cc}=\mathcal{G}$. By [C] Corollary 4.4.4.12 a), $\f^{cc}$ is the W*-subalgebra of $\lb{E}{H}$ generated by $\f$.\qed
 
\begin{center}
{\bfseries REFERENCES}
\end{center}
\begin{flushleft}
[C] Corneliu Constantinescu, {\it C*-algebras.} Elsevir, 2001. \newline
\end{flushleft}
\begin{flushright}
{\scriptsize \hspace{-5mm} Corneliu Constantinescu$\quad$\\
Bodenacherstr. 53$\qquad\;$\\
CH 8121 Benglen$\qquad\;\;$\\
e-mail: constant@math.ethz.ch }
\end{flushright}
\end{document}